\documentclass[]{amsart}
\usepackage{amssymb,mathdots}
\usepackage{mathrsfs}
\usepackage[applemac]{inputenc}
\usepackage{amsmath, graphicx}
\usepackage{wrapfig}
\usepackage{enumerate}
\usepackage{url}

\def\scaledpicture#1by#2(#3scaled#4){{
\dimen0=#1  \dimen1=#2
\divide\dimen0 by 1000 \multiply\dimen0 by #4
\divide\dimen1 by 1000 \multiply\dimen1 by #4
\picture \dimen0 by \dimen1 (#3 scaled #4)}}
\def\dfigure#1by#2(#3scaled#4offset#5:#6)
    {\medskip
     \vglue 2mm minus 2mm
     $$
       \hbox{
         \hglue#5
         {\scaledpicture #1 by #2 (#3 scaled #4)}
       }
     $$
     \par\goodbreak
     \vglue 2mm minus 2mm
     \medskip}

\def\qmod#1#2{{\hbox{}^{\displaystyle{#1}}}\!\big/\!\hbox{}_{
\displaystyle{#2}}}

\def\resto#1#2{{
#1\hskip 0.4ex\vline_{\hskip 0.2ex\raisebox{-0,2ex}
{{${\scriptstyle #2}$}}}}}

\def\C{{\mathbb C}}

\def\P{{\mathbb P}}

\def\R{{\mathbb R}}
\def\Z{{\mathbb Z}}

\def\map{\longrightarrow}
\def\textmap#1{\mathop{\vbox{\ialign{
                                  ##\crcr
      ${\scriptstyle\hfil\;\;#1\;\;\hfil}$\crcr
      \noalign{\kern 1pt\nointerlineskip}
      \rightarrowfill\crcr}}\;}}
\def\bigtextmap#1{\mathop{\vbox{\ialign{
                                  ##\crcr
      ${\hfil\;\;#1\;\;\hfil}$\crcr
      \noalign{\kern 1pt\nointerlineskip}
      \rightarrowfill\crcr}}\;}}
      
\newcommand{\cal}{\mathcal}
\def\textlmap#1{\mathop{\vbox{\ialign{
                                  ##\crcr
      ${\scriptstyle\hfil\;\;#1\;\;\hfil}$\crcr
      \noalign{\kern-1pt\nointerlineskip}
      \leftarrowfill\crcr}}\;}}

\def\g{{\mathfrak g}}

\newtheorem{sz}{Satz}%[section]
\newtheorem{thry}[sz]{Theorem}

\newtheorem{re}[sz]{Remark}
\newtheorem{co}[sz]{Corollary}
\newtheorem{dt}[sz]{Definition}
\newtheorem{lm}[sz]{Lemma}

\begin{document}
 
\def\tr{\mathrm {Tr}}
\def\End{\mathrm {End}}
\def\Aut{\mathrm {Aut}}
\def\Spin{\mathrm {Spin}}
\def\U{\mathrm{U}}
\def\SU{\mathrm {SU}}
\def\SO{\mathrm {SO}}
\def\PU{\mathrm {PU}}
\def\GL{\mathrm {GL}}
\def\spin{\mathrm {spin}}
\def\u{\mathrm {u}}
\def\su{\mathrm {su}}
\def\so{\mathrm {so}}
\def\pu{\mathrm {pu}}
\def\Pic{\mathrm {Pic}}
\def\Iso{\mathrm {Iso}}
\def\NS{\mathrm{NS}}
\def\deg{\mathrm {deg}}
\def\Hom{\mathrm{Hom}}
\def\Herm{\mathrm{Herm}}
\def\Vol{{\rm Vol}}
\def\pf{{\bf Proof: }}
\def\id{ \mathrm{id}}
\def\Im{\mathrm{Im}}
\def\im{\mathrm{im}}
\def\rk{\mathrm {rk}}
\def\ad{\mathrm {ad}}
%\def\coker{\mathrm{coker}}
%%%%%%%%%%%%%%%%%%%%%%%%%%
\def\spc{\mathrm{Spin}^c}
\def\U2{\mathrm{U(2)}}
\def\niq{=\kern-.18cm /\kern.08cm}
%%%%%%%%%%%%%%%%%%%
\def\Ad{\mathrm {Ad}}
\def\RSU{\R\mathrm{SU}}
\def\ad{{\rm ad}}
\def\dva{\bar\partial_A}
\def\da{\partial_A}
\def\p{{\rm p}}
\def\sp{\Sigma^{+}}
\def\sm{\Sigma^{-}}
\def\spm{\Sigma^{\pm}}
\def\smp{\Sigma^{\mp}}
\def\oo{{\scriptstyle{\cal O}}}
\def\ooo{{\scriptscriptstyle{\cal O}}}
\def\sw{Seiberg-Witten }
\def\pa{\partial_A\bar\partial_A}
\def\Dr{{\raisebox{0.17ex}{$\not$}}{\hskip -1pt {D}}}
\def\gr{{\scriptscriptstyle|}\hskip -4pt{\g}}
\def\subsetint{{\  {\subset}\hskip -2.45mm{\raisebox{.28ex}
{$\scriptscriptstyle\subset$}}\ }}

\title
{Symmetric theta divisors of Klein surfaces}
\author{Christian Okonek \& Andrei Teleman}
 \begin{abstract}  This is a slightly expanded version of the talk given by Ch.O. at the conference ``{\it Instantons in complex geometry}", at the Steklov Institute in Moscow. The purpose of this talk was to explain the algebraic results of our paper "{\it Abelian Yang-Mills theory on Real tori and Theta divisors of Klein surfaces}" \cite{OT2}.

In this paper we compute  determinant index bundles of certain families of Real Dirac type operators on Klein surfaces as elements in the corresponding Gro\-then\-dieck  group \cite{G} of Real line bundles in the sense of Atiyah. On a Klein surface these determinant index bundles have a natural holomorphic description  as theta line bundles. In particular we compute the first Stiefel-Whitney classes of the corresponding fixed point bundles on the real part 
 of the Picard torus. The
computation of these classes is important, because they control to a large
extent the orientability of certain moduli spaces in Real gauge theory and Real algebraic geometry \cite{OT3}.
 
 \end{abstract} 
\maketitle

Let $C$ be a Riemann surface, and $g:=h^0(\omega_C)$ its genus. {\it The   geometric theta di\-vi\-sor} of $C$ is the effective divisor of $\Pic^{g-1}(C)$ defined by
$$\Theta:=\{[\cal L]\in\Pic^{g-1}(C)|\ h^0(\cal L)>0\}\ .
$$

The set $\theta$ of {\it theta characteristics} of $C$ is the set of square roots of $[\omega_C]$, i.e.
$$\theta:=\{[\kappa]\in\Pic^{g-1}(C)|\ \kappa^{\otimes 2}\simeq\omega_C\}\subset\Pic^{g-1}(C)\ .
$$
Note that $\# \theta=2^{2g}$. For a theta characteristic $[\kappa]\in\theta$ we obtain an associated {\it symmetric theta divisor}
$$\Theta_{[\kappa]}:=\Theta-[\kappa]\subset\Pic^0(C)
$$
and a holomorphic line bundle ${\cal L}_{[\kappa]}:={\cal O}_{\Pic^0(C)}(\Theta_{[\kappa]})$ on the Abelian variety $\Pic^0(C)$.
\vspace{2mm}
\begin{dt} A Klein surface is a pair $(C,\iota)$, where $C$ is a Riemann surface and $\iota:C\to C$ an anti-holomorphic involution. 
\end{dt} 
\vspace{2mm}
The fixed point locus $C^\iota$ decomposes is a finite union $\coprod_{i=1}^n C_i$ of $n:=\#\pi_0(C^\iota)$ circles $C_i$. The topological type of the Klein surface $(C,\iota)$ is the triple  $(g,n,a)$, where
$$
a:=\left\{\begin{array}{cccc}
0&\rm when& \qmod{C}{\langle\iota\rangle}&\hbox{is orientable}\\
1&\rm when& \qmod{C}{\langle\iota\rangle}&\hbox{is non-orientable.}
\end{array}\right.
$$
\\

These invariants are subject to the following conditions \cite{GH}:
\begin{enumerate}
\item $1\leq n\leq g+1$, $g+1-n\equiv 0$ (mod 2), when $a=0$,
\item  $0\leq n\leq g$, when $a=0$.
\end{enumerate}
\vspace{3mm}
{\bf Examples:}
\begin{enumerate}
\item $g=0$, so $C=\P^1_\C$.
\begin{enumerate}
\item   $a=0$: $\iota[z_0,z_1]=[\bar z_0,\bar z_1]$; in this case one has $C^\iota=\P^1_\R$, and $C/\langle\iota\rangle\simeq D^2$,
\item  $a=1$: $\iota[z_0,z_1]=[-\bar z_1,\bar z_0]$; in this case one has $C^\iota=\emptyset$, $C/\langle\iota\rangle\simeq \P^2_\R$
\end{enumerate}
\item $g=1$, so $C$ is an elliptic curve, say $C=\C/\langle 1,\tau\rangle$ with $\tau\in\C\setminus\R$.
\begin{enumerate}
\item   $a=0$: take $\tau=it$ with $t\in\R_{>0}$, $\iota[z]=[\bar z]$; in this case $C^\iota=C_1\coprod C_2$, and $C/\langle\iota\rangle$  is an annulus.
\item   $a=1$: take $\tau=it$ with $t\in\R_{>0}$, $\iota([z])=[\frac{1}{2}+\bar z]$; in this case $C^\iota=\emptyset$, and $C/\langle\iota\rangle$  is a Klein bottle.
\item   $a=1$: take $\tau=\frac{1}{2}+it$ with $t\in\R_{>0}$, $\iota([z])=[\bar z]$; in this case $C^\iota=C_1$, and $C/\langle\iota\rangle$  is a Möbius band.
\end{enumerate}

%\begin{wrapfigure}[8]{r}{.5\textwidth}
%\vspace*{-4mm}
%\hspace*{0.6cm}
%\includegraphics[scale=0.3]{riemann.pdf}
%\caption{Legende de la figure}%\label{01_img_01}
%\end{wrapfigure}

\item The picture below represents a Klein surface with topological type $(4,3,0)$
\vspace{2mm}
\begin{center} 
\includegraphics[width=6cm]{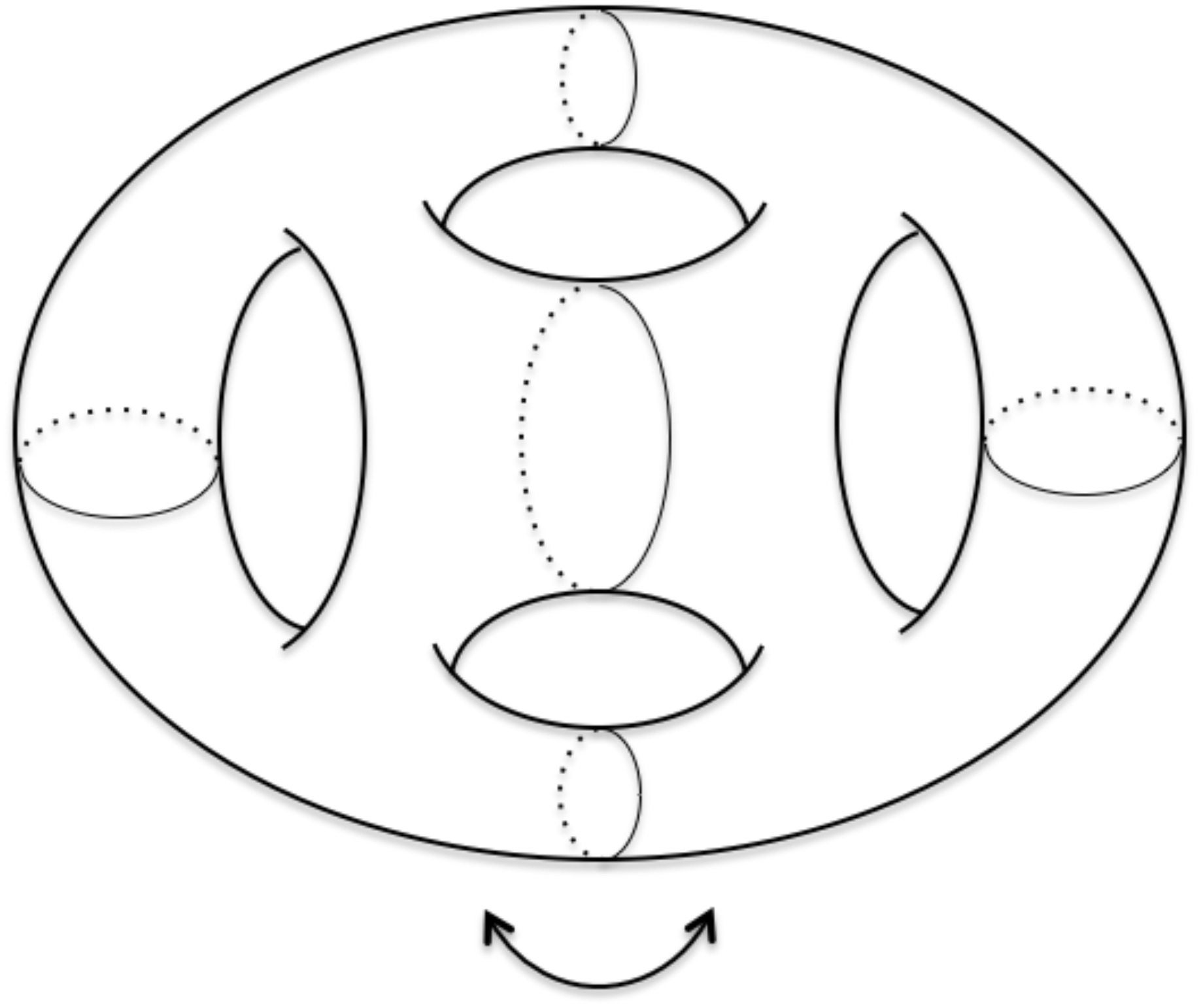}
\end{center}
\end{enumerate}
\vspace{2mm}
\begin{re} Let $X$ be a compact, connected complex manifold and $\tau:X\to X$ an anti-holomorphic involution. Then $\tau$ induces an anti-holomorphic involution
$$\hat \tau:\Pic(X)\to\Pic(X)
$$
defined by $\hat \tau([\cal L]):=[\tau^*(\bar{\cal L})]$.
\end{re}
\vspace{2mm}

Let $(X,\tau)$ be a  compact, connected  complex manifold endowed with an anti-holomorphic involution such that $X^\tau\ne\emptyset$. We define a morphism $w:\Pic(X)^{\hat\tau}\to H^1(X^\tau,\Z_2)$ in the following way:  for a holomorphic line bundle ${\cal L}$ with $\tau^*(\bar{\cal L})\simeq\cal L$ there exists an   involutive, anti-holomorphic lift $\tilde\tau_{\cal L}:{\cal L}\to{\cal L}$ of $\tau$, which is unique up to multiplication by constants $\zeta\in S^1$. The fixed point locus  ${\cal L}^{\tilde\tau_{\cal L}}$ can be regarded as a real line bunde over $X^\tau$, and we put
$$w([\cal L]):=w_1({\cal L}^{\tilde\tau_{\cal L}})\ .
$$

We are interested in the following
\\ \vspace{1mm}\\
{\bf Problem}: For a Klein surface $(C,\iota)$ with $C^\iota\ne\emptyset$, let $[\kappa]\in\theta$ such that $\hat\iota[\kappa]=[\kappa]$, and let ${\cal L}_{[\kappa]}$ be the associated holomorphic line bundle on $\Pic^0(C)$. Compute
$$w({\cal L}_{[\kappa]})\in H^1(\Pic^0(C)^{\hat\iota},\Z_2)\ .
$$
\vspace{2mm}\\

Our problem is motivated by Real gauge theory, namely by Real Gromov-Witten and Real Seiberg-Witten theory. The point is that ${\cal L}_{[\kappa]}$ can   be regarded  as the determinant line bundle associated with a family of order 0 perturbations of the  Dirac operator associated with the Spin structure on $C$ defined by $[\kappa]$. This family of perturbations is parameterized by $\Pic^0(C)$. Using this interpretaion, one sees that $w({\cal L}_{[\kappa]})$ controls the orientability of different components of Real moduli spaces of generalized vortices.
\vspace{2mm}
\begin{re} The Chern class $c_1({\cal L}_{[\kappa]})$ can be computed  using the Grothendieck-Riemann-Roch theorem applied to the projection $\Pic^0(C)\times C\to\Pic^0(C)$ and the interpretation of ${\cal O}_{\Pic^{g-1}(C)}(\Theta)$ as the determinant line bundle of the total direct image of a Poincaré line bundle on $\Pic^{g-1}(C)\times C$.
\end{re}
\vspace{2mm}

In order to solve our problem we shall use the following strategy: we  show that $w({\cal L}_{[\kappa]})$ can be read from the Appell-Humbert data of ${\cal L}_{[\kappa]}$ regarded as holomorphic line bundle on the torus $H^1(C,{\cal O}_C)/H^1(X,\Z)$. Therefore, our strategy has two steps:
\begin{enumerate}
\item[I)]\label{1step} Determine explicitly the Appell-Humbert data of ${\cal L}_{[\kappa]}$,\vspace{2mm}
\item[II)] \label{2step} Extract the Stiefel-Whitney class $w({\cal L}_{[\kappa]})$ from these data.
\end{enumerate}
\vspace{4mm}

The general set-up is the following:  $V$ is a finite  dimensional complex vector space, $\Lambda\subset\C$ a maximal lattice, and $\tau:V\to V$ an anti-linear involution such that $\tau(\Lambda)\subset \Lambda$. The holomorphic torus $T:=V/\Lambda$ comes with an induced involution which will be denoted by the same symbol $\tau$.
\vspace{2mm}
\begin{dt} An Appel-Humbert datum for the pair $(V,\Lambda)$ is a pair $(H,\alpha)$, where $H:V\times V\to\C$ is a Hermitian form on $V$  such that $(\im H)(\Lambda\times\Lambda)\subset\Z$, and $\alpha:\Lambda\to S^1$ is an $\im H$--semi-character, i.e. it satisfies the identity
$$\alpha(\lambda+\lambda')=\alpha(\lambda)\alpha(\lambda') e^{\pi i\ \im H(\lambda,\lambda')}  \ . 
$$
\end{dt}
\vspace{2mm}

We  denote by $\Hom_{\im H}(\Lambda,S^1)$ the set of $\im H$--semi-characters on $\Lambda$. An Appel-Humbert datum $(H,\alpha)$ defines {\it a canonical factor of automorphy} 
$$a(H,\alpha):\Lambda\times V\to\C^*$$ defined by
$$a(H,\alpha)(\lambda,v):=\alpha(\lambda) e^{\pi[H(\lambda,v)+\frac{1}{2} H(\lambda,\lambda)] }\ .
$$
Let ${\cal L}(H,\alpha)$ be the holomorphic line bundle   ${\cal L}(H,\alpha):=V\times\C/\sim_{a(H,\alpha)}$, where $\sim_{a(H,\alpha)}$ is the equivalence relation on the product $V\times\C$ defined by the factor of automorphy $a(H,\alpha)$.
\vspace{2mm}
\begin{thry} (Appel-Humbert) The assignement $(H,\alpha)\mapsto {\cal L}(H,\alpha)$ induces an isomorphism
$$AH: \coprod_{\begin{array}{c}  
\scriptstyle{H\ Hermitian\  form\ on\ V}\vspace{-1mm}\\
\scriptstyle(\im H)(\Lambda\times\Lambda)\subset\Z\end{array}}
\hspace{-12mm}\Hom_{\im H}(\Lambda,S^1)\to\Pic(T)\ .
$$
\end{thry}
\vspace{2mm}

We give a new proof of this statement which will allow us to give an important, new geometric interpretation of the Appel-Humbert datum corresponding to a given holomorphic line bundle.
This new proof is based on two ideas: \vspace{2mm}
\begin{enumerate}
\item[1)] Prove a  general result about the holonomy of Yang-Mills connections on line bundles,
\item[2)] Use the Kobayashi-Hitchin correspondence to recover the Appel-Humbert theorem.
\end{enumerate} 
\vspace{2mm}
1) Let $u\in \mathrm{Alt}^2(\Lambda,\Z)=H^2(T,\Z)$ and $L_u$ be a Hermitian line bundle on $T$ with $c_1(L_u)=u$. The linear extension $u_\R:V\times V\to\R$  of $u$ can also be regarded as a differentiable 2-form on $T$, which is harmonic with respect to the  flat metric induced by an inner product on $V$. The moduli space of Yang-Mills connections on $L_u$ is
$${\cal M}(L_u):=\left\{[A]\in\qmod{{\cal A}(L_u)}{\cal G}\ \vline\ \ \frac{i}{2\pi} F_A=u_\R\right\},\ 
$$
where ${\cal A}(L_u)$ stands for  space of unitary connections on $L_u$, and  ${\cal G}:={\cal C}^\infty(T,S^1)$ denotes the gauge group of the Hermitian line bundle $L_u$.

A unitary connection $A\in{\cal A}(L_u)$ defines a map $\alpha^A:\Lambda\to S^1$ given by the holonomy of $A$ along the loops $c_\lambda$ in $T$ corresponding to the  lattice elements $\lambda\in\Lambda$. More precisely 
$$\alpha^A(\lambda):=h^A_{c_\lambda}\in S^1\ ,
$$
where $c_\lambda:S^1=\R/\Z\to T$ is defined  by $c_\lambda(t):=[t\lambda]$.
\vspace{2mm}
\begin{thry} The assignment $A\mapsto\bar\alpha^A$ induces a homeomorphism 
$$h:{\cal M}(L_u)\to\Hom_u(\Lambda,S^1)\ .$$
\end{thry}
\vspace{1mm}{\ }\\
2) Let $u\in\mathrm{Alt}^2(\Lambda,\Z)$ such that the associated real form $u_\R$ has type $(1,1)$, i.e. $J^*(u_\R)=u_\R$. Denote by $H_u$ the associated Hermitian form  
$$H_u(v,w):=u_\R(v,Jw)+i u_\R(v,w) \ .
$$
We get a commutative diagram:\\ \\
$$
\begin{array}{ccc}\hspace{8mm}
\displaystyle\coprod_{\begin{array}{c}\scriptstyle H\ { Hermitian\  form\ on\ }V\vspace{-1mm}\\\scriptstyle (\im H)(\Lambda\times\Lambda)\subset\Z\end{array}}\hspace{-12mm}\Hom_{\im H}(\Lambda,S^1)&=&
\displaystyle\coprod_{\begin{array}{c}
\scriptstyle u\in\mathrm{Alt}^2(\Lambda,\Z)\vspace{-1mm}\\\scriptstyle J^*(u_\R)=u_\R\end{array}}\hspace{-6mm} \Hom_{\im H_u}(\Lambda,S^1)\\
\\
{\ }\downarrow AH&& \downarrow h^{-1}\ \ \ \ \ \ \ \ \ \ \ \\
\\
\Pic(T)&\stackrel{KH}{\longleftarrow}\ &\displaystyle\coprod_{\scriptstyle u\in \mathrm{NS}(T)^{\phantom{I^I}}} \hspace{-4mm} {\cal M}(L_u)\ \ \ \ \ \ \ \  \ \  
\end{array}
$$
\vspace{1mm}\\

Here $KH$ denotes the Kobayashi-Hitchin correspondence between equivalence classes of Hermite-Einstein connections and isomorphism classes of (polystable) holomorphic vector bundles. Since any holomorphic line bundle is stable, the polystability condition is empty in our situation. Recall that, in general,  a unitary connection on a Hermitian vector bundle on a compact Kähler manifold is Hermite-Einstein if and only if it is Yang-Mills and its curvature has type (1,1).
\vspace{2mm}
\begin{re} The commutative diagram above gives a geometric interpretation of the AH data $(H,\alpha)$ of a holomorhic line bundle ${\cal L}$: $(H,\alpha)$ corresponds to the curvature and the holonomy of the (essentially unique) Hermite-Einstein connection which is compatible with the holomorphic structure of  ${\cal L}$.
\end{re} 
\vspace{2Mm}

Now we come  back to our problem: for a $\hat\tau$-invariant element $[\cal L]\in\Pic(T)^{\hat\tau}$ we want to determine $w({\cal L})\in H^1(T^\tau,\Z_2)$. Note that the fixed point locus $T^\tau$ decomposes as a disjoint union
$$T^\tau=\coprod_{\raisebox{-1mm}{$\scriptstyle[\mu]\in$ }\frac{1}{2}\Lambda^{-\tau} /\hspace{-0.5mm}\raisebox{-2mm}{$
\scriptstyle\frac{1}{2}(1-\tau)\Lambda$ }}\ \left(\qmod{V^\tau}{\Lambda^\tau}+[\mu]\right)\ ,
$$
so that $w(\cal L)$ can be regarded as a map
$$w(\cal L)\ :\ \qmod{\frac{1}{2}\Lambda^{-\tau}}{\frac{1}{2}(1-\tau) \Lambda}\map \Hom(\Lambda^\tau,\Z_2)\ .
$$
Using this point of view one can prove:
\vspace{2mm}
\begin{thry} Let $(H,\alpha)$ be an AH datum, and put ${\cal L}:={\cal L}(H,\alpha)$.  When $[{\cal L}]$ is $\hat\tau$-invariant,  one has:
\begin{enumerate}
\item $w({\cal L})([\mu])= w({\cal L})([0])+\overline{\im H(2\mu,\cdot)}$ for every  $\mu\in\frac{1}{2}\Lambda^{-\tau}$,
\item $w({\cal L})([\mu])(\lambda+\tau\lambda)=\overline{\im H(\lambda,\tau\lambda)}$ for every $\lambda\in\Lambda$,
\item $w(\cal L)([0])(\lambda)=\bar\alpha(\lambda)$ for every $\lambda\in\Lambda^\tau$.
\end{enumerate}
\end{thry}
\vspace{2mm}
In this statement the bar on the right means congruence class (mod 2) in the first two formulae and conjugation in the third. This result completes step  II) of our strategy, namely it allows us to extract the map $w$ associated to a fixed point of $\hat\tau$ from the corresponding Appel-Humbert datum. We now have  to  complete step I), namely to find the Appel-Humbert datum of the line bundles ${\cal L}_{[\kappa]}$  associated with an $\hat\iota$-invariant theta characteristic $[\kappa]$.

Recall that, according to Mumford \cite{ACGH}, a theta characteristic $[\kappa]\in\theta$ has an associated   theta form $q_{[\kappa]}:\Pic^0(C)_2\to\Z_2$ defined on the 2-torsion subgroup $\Pic^0(C)_2$ of $\Pic^0(C)$  and given by
$$q_{[\kappa]}([\eta]):= h^0(\kappa\otimes\eta)-h^0(\kappa)\ ({\rm mod}\ 2)\ .
$$
Identifying $\Pic^0(C)_2\simeq H_1(C,\Z_2)$ this form $q_{[\kappa]}$  satisfies the {\it Riemann-Mumford} relations:
\begin{equation}\label{RM}
q_{[\kappa]}([\eta]+[\eta'])=q_{[\kappa]}([\eta])+q_{[\kappa]}([\eta'])+[\eta]\cdot[\eta']\ .
\end{equation}
\vspace{-1mm}
\begin{thry} Let $[\kappa]\in\theta$. Then ${\cal L}_{[\kappa]}\simeq{\cal L}(H_{\langle\cdot,\cdot\rangle},\alpha_{[\kappa]})$, where 
\begin{enumerate}
\item[i)] $\langle\cdot,\cdot\rangle:H^1(C,\Z)\times H^1(C,\Z)\to\Z$ is the cup form,
\item[ii)] $\alpha_{[\kappa]}$ is defined by the formula $\alpha_{[\kappa]}(\lambda):=(-1)^{q_{[\kappa]} (\overline{\lambda\cap[C]} )}$.
\end{enumerate}
\end{thry}{\ }\vspace{-3mm}\\
{\it Idea of proof:}
Since $c_1({\cal L}_{[\kappa]})=\langle\cdot,\cdot\rangle$ it follows that the first component of the AH datum of ${\cal L}_{[\kappa]}$ is $H_{\langle\cdot,\cdot\rangle}$ as claimed. On the other hand, using the fact that the divisor $\Theta_{[\kappa]}\subset\Pic^0(C)$ is symmetric (i.e. $(-1)^*\Theta_{[\kappa]}=\Theta_{[\kappa]}$), one knows \cite{BL}, that the second component of the AH datum is the $\langle\cdot,\cdot\rangle$-semicharacter $\alpha_{[\kappa]}$ given by
$$\alpha_{[\kappa]}=(-1)^{{\rm mult}_{[\frac{1}{2}\lambda]}(\Theta_{[\kappa]})-{\rm mult}_{[0]}(\Theta_{[\kappa]})}\ .
$$
Now we use Riemann's singularity theorem \cite{BL}, which states 
$${\rm mult}_{[{\cal L}]}\Theta=h^0(\cal L)\ ,
$$
and the identification $\frac{1}{2}H^1(C,\Z)/H^1(C,\Z)\simeq H^1(C,\Z_2)$ given by 
$$[\frac{1}{2}\lambda]\mapsto \overline{\lambda\cap[C]}\ .$$
\vspace{-1mm}
\begin{co} Let $[\kappa]$ be an $\hat\iota$-invariant theta characteristic.  Then
$$w({\cal L}_{[\kappa]})([0])(\lambda)=(-1)^{q_{[\kappa]}(\overline{\lambda\cap[C]})}\ \forall\lambda\in H^1(C,\Z)^{-\iota^*}\ .
$$
\end{co}
\vspace{2mm}

This result is not entirely satisfactory, because the right hand side is not topological. We need a third, a priori unexpected  step, which will give  an explicit,  purely topological interpretation of the right hand term.
First of all recall  that 
\begin{equation}\label{trivial}
w({\cal L}_{[\kappa]})([0])(\lambda-\iota^*(\lambda))=\overline{\langle\lambda,-\iota^*\lambda\rangle}\ \forall\lambda\in H^1(C,\Z)\ .\end{equation}
Using results of   \cite{CN} the following is easy to see:
\vspace{2mm}
\begin{lm} \label{gen} Let $C^\iota=\coprod_{i=1}^n C_i$ be the decomposition of the fixed point locus of $\iota$ in connected components,  choose orientations of these components, and denote by $[C_1]^\vee$ the cohomology classes which correspond  to $[C_i]$ via Poincaré duality. Then   
$$\langle [C_1]^\vee,\dots,[C_n]^\vee\rangle_{\Z}\ ,\ (1-\iota^*) H^1(C,\Z)
$$
generate $H^1(C,\Z)^{-\iota^*}$.
\end{lm}
\vspace{2mm}
Combining (\ref{trivial}) and Lemma \ref{gen} we see that it suffices to compute
$$q_{[\kappa]}(\overline{[C_i]^\vee\cap[C]})=q_{[\kappa]}([C_i]_2)\ .
$$
\vspace{-1mm}
\begin{thry} Let $(C,\iota)$ be a Klein surface with $C^\iota=\coprod_{i=1}^n C_i$, where $n>0$. Let $[\kappa]$ be an $\hat\iota$-invariant theta characteristic. Then
$$w_1(\resto{{\cal L}_{[\kappa]}^{\tilde\iota_{{\cal L}_{[\kappa]}}})}{\Pic^0(C)^{\hat\iota}_0})([C_i]^\vee)=(-1)^{\langle w_1(\kappa^{\tilde\iota_\kappa}),[C_i]_2 \rangle+1}\ .
$$
\end{thry}
\vspace{2mm}
Here we have denoted by $\Pic^0(C)^{\hat\iota}_0$ the connected component of  the trivial line bundle $[{\cal O}_C]$ in the fixed point locus $\Pic^0(C)^{\hat\iota}$.
{\ }\\ \\ {\it Idea of proof:} We have to show that
$$q_{[\kappa]}([C_i]_2)=\langle w_1(\kappa^{\tilde\iota_\kappa}),[C_i]_2\rangle+1\ .
$$
which gives a topological interpretation of Mumford's theta form.  The proof of this formula uses delicate topological arguments combined with results of Johnson \cite{J}, Libgober \cite{L},  and Atiyah \cite{A2}. To explain their results, consider the diagram:
$$
\begin{array}{rcl}\vspace{2mm}
\theta & \hspace{-6mm} \stackrel{\xi}{\map}&\hspace{-10mm} \mathrm{Spin}(C)\\
&\hspace{-2mm}\raisebox{-1mm}{$q$}\hspace{-1mm}\searrow \hspace{5mm}\ \ \ \  \ \swarrow\omega
\vspace{2mm}\\
&\hspace{-2mm} Q(H_1(C,\Z_2),\cdot) 
\end{array}  
$$
Here $Q(H_1(C,\Z_2),\cdot)$ denotes the set of maps $q:\Pic^0(C)_2\to\Z_2$  satisfying the Riemann-Mumford relations (\ref{RM}), $\mathrm{Spin}(C)$ is the set of equivalence classes of $\Spin$-structures on $C$, $\xi$ is the  correspondence between theta charcateristics and $\Spin(C)$ defined by Atiyah, and $\omega$ is a map defined by Johnson \cite{J} in purely topological terms. The result follows from the fact that $\xi$ is bijective \cite{A2}, $\omega$ is bijective \cite{J}, and $\omega\circ\xi=q$ \cite{L}.

{\ }
\vspace{10mm}  \\
{\small Christian Okonek: \\
Institut f\"ur Mathematik, Universit\"at Z\"urich,
Winterthurerstrasse 150, CH-8057 Z\"urich,\\
e-mail: okonek@math.unizh.ch
\\  \\
Andrei Teleman: \\
CMI,   Universit\'e de Provence,  39  Rue F. Joliot-Curie, F-13453
Marseille Cedex 13,   e-mail: teleman@cmi.univ-mrs.fr
}

  \end{document}